\documentclass{proc-l}

\newtheorem{theorem}{Theorem}[section]
\newtheorem{lemma}[theorem]{Lemma}

\theoremstyle{definition}
\newtheorem{definition}[theorem]{Definition}

\theoremstyle{remark}
\newtheorem{remark}[theorem]{Remark}
\newtheorem*{acknowledgements}{Acknowledgements}
\newtheorem{corollary}[theorem]{Corollary}
\numberwithin{equation}{section}

\def\openone

{\mathchoice

{\hbox{\upshape \small1\kern-3.3pt\normalsize1}}

{\hbox{\upshape \small1\kern-3.3pt\normalsize1}}

{\hbox{\upshape \tiny1\kern-2.3pt\SMALL1}}

{\hbox{\upshape \Tiny1\kern-2pt\tiny1}}}

\makeatletter

\newbox\ipbox

\newcommand{\ip}[2]{\left\langle #1\,|\,#2\right\rangle}

\newcommand{\diracb}[1]{\left\langle #1\mathrel{\mathchoice

{\setbox\ipbox=\hbox{$\displaystyle \left\langle\mathstrut 
#1\right.$}

\vrule height\ht\ipbox width0.25pt depth\dp\ipbox}

{\setbox\ipbox=\hbox{$\textstyle \left\langle\mathstrut 
#1\right.$}

\vrule height\ht\ipbox width0.25pt depth\dp\ipbox}

{\setbox\ipbox=\hbox{$\scriptstyle \left\langle\mathstrut 
#1\right.$}

\vrule height\ht\ipbox width0.25pt depth\dp\ipbox}

{\setbox\ipbox=\hbox{$\scriptscriptstyle \left\langle\mathstrut 
#1\right.$}

\vrule height\ht\ipbox width0.25pt depth\dp\ipbox}

}\right. }

\newcommand{\dirack}[1]{\left. \mathrel{\mathchoice

{\setbox\ipbox=\hbox{$\displaystyle \left.\mathstrut 
#1\right\rangle$}

\vrule height\ht\ipbox width0.25pt depth\dp\ipbox}

{\setbox\ipbox=\hbox{$\textstyle \left.\mathstrut 
#1\right\rangle$}

\vrule height\ht\ipbox width0.25pt depth\dp\ipbox}

{\setbox\ipbox=\hbox{$\scriptstyle \left.\mathstrut 
#1\right\rangle$}

\vrule height\ht\ipbox width0.25pt depth\dp\ipbox}

{\setbox\ipbox=\hbox{$\scriptscriptstyle \left.\mathstrut 
#1\right\rangle$}

\vrule height\ht\ipbox width0.25pt depth\dp\ipbox}

} #1\right\rangle}

\newcommand{\ltwor}{L^{2}\left(\mathbb{R}\right)}

\newcommand{\cj}[1]{\overline{#1}}

\newcommand{\bz}{\mathbb{Z}}
\newcommand{\br}{\mathbb{R}}

\usepackage{graphicx}

\newcommand{\D}{\mathbf{D}}
\newcommand{\T}{\mathbf{T}}
\newcommand{\ltword}{L^2(\mathbb{R}^d)}
\newcommand{\U}{\mathbf{U}_P}
\begin{document}

\title{Oversampling generates super-wavelets}


\author{Dorin Ervin Dutkay}
\address{ Department of Mathematics,
Hill Center-Busch Campus,
Rutgers, The State University of New Jersey,
110 Frelinghuysen Rd,
Piscataway, NJ 08854-8019, USA}
\curraddr{}
\email{ddutkay@math.rutgers.edu}
\thanks{}

\author{Palle Jorgensen}
\address{Department of Mathematics, The University of Iowa, 14 MacLean Hall, Iowa City, IA, 52242, USA}
\curraddr{}
\email{jorgen@math.uiowa.edu}
\thanks{}

\subjclass[2000]{42C40, 47A20, 65T60, 94A20}
\keywords{wavelet, frame, sampling, oversampling, affine, scaling, lattice, interpolation, dilations, extensions, super wavelets, operators, frames, Hilbert space}

\date{November 15, 2005}

\dedicatory{}

\commby{}

\begin{abstract}
We show that the second oversampling theorem for affine systems generates super-wavelets. These are frames generated by an affine structure on the space $\underbrace{\ltword\oplus...\oplus\ltword}_{p\mbox{ times}}$.
\end{abstract}

\maketitle


\newcommand{\tr}{\theta_r}
\newcommand{\tts}{\theta^*}
\newcommand{\ttq}{\theta^*_q}
\newcommand{\sg}{\sigma}
\newcommand{\sgs}{\sigma^*}

\section{Introduction}
\par
While, as is well known the study of wavelets draws from a variety of areas of mathematics, and wavelet algorithms have numerous applications, in this paper we focus on an operator theoretic aspect of the subject. We feel that our approach clarifies fundamental techniques in the subject, and at the same time may be of independent interest in operator theory and in analysis. Some key ideas begin with early engineering applications; e.g., to speech, and to time-series. Here ``oversampling'' refers to certain redundancies that are introduced into signal processing. When our problem is formulated in the context of a fixed Hilbert space, we have a useful notion of dilation (or extension) available from operator theory; that of passing to a bigger Hilbert space $H$ (a ``super space'') where bases and discrete transforms work without redundancies, and where as a result computations simplify. Even certain symmetries are better understood in the super space $H$. In the end, answers can be restricted back to smaller Hilbert space again.
\par
The subject of wavelets originates with a dual track of algorithms; one in function theory, and one in signal processing. Examples: subband coding, pyramid algorithms, multiresolutions, adaptive bases, and data compression with thresholding. 
\par
First recall from the wavelet literature (e.g., \cite{Du})  that fundamentals from signal processing may be thought of and understood within the framework of operator algebra and representation theory. In fact, this viewpoint is almost ubiquitous and it may be used to advantage in a wider context of mathematical analysis ; i.e., used in processes which select and analyze special bases in function spaced. It applies even more generally when the notion of an orthonormal basis (ONB) in Hilbert space is extended in such a way to allow instead {\it frames}. A frame system is defined rigorously in equation (\ref{eqfr}) below. Within the family of wavelet bases, frames have the advantage of encompassing wavelet bases that allow more symmetries that can be accommodated by ONB-wavelets. This is critical to applications, for example to the algorithms used in digital fingerprint compression, see e.g., \cite{Bri95}. For these, symmetry is a critical feature, and so orthogonality must be relaxed. (See e.g., \cite[Chapter 6]{Dau92}.) The bases we consider in this paper serve as framework for both discrete signal processing and more function theoretic models. A leading theme in our analysis is the interplay between the two. A key link between the two is provided by time-frequency analysis. We begin in section 3 with a formulation of Fourier duality which is tailored to our presentation of ``oversampling'' as part of super-wavelets. The term ``oversampling'' is from signal processing and coding theory. Starting with a prescribed sample point, there are reasons for ``oversampling'', i.e., the addition of more sample points, or generation of redundancies: for example, the redundancies may allow error detection and/or error correction. A concrete way to pass to a desired oversampling is to start with sample points which are restricted to a fixed lattice $L$ in $\br^d$, i.e., a rank-$d$ abelian and discrete subgroup of $\br^d$. It follows that a bigger lattice $L'$ of sample points may then be represented in the form   $P L' = L$ where $P$ is an invertible and integral $d$ by $d$ matrix; i.e., $L' = P^{-1} L$. In the case a frame system is first created in the Hilbert space $L^2(\br^d)$ from a suitably chosen pair $(L, P)$ where $L$ is a lattice and $P$ is a scaling matrix, we show that the super-Hilbert space $H$ may then be realized as a $p$-fold orthogonal sum of $\ltword$ with itself where $p = |\det P|$. Moreover we give explicit formulas for the respective actions in the two Hilbert spaces which in turn are expressed and explained in terms of ``oversampling''.

\par
We use the term sampling in the sense of Shannon. Shannon showed that a function $f$ on the real line $\br$ which has its Fourier transform $\hat f$ supported in a bounded interval can be reconstructed by interpolation of its values resulting from sampling at integral multiples of a certain rate $\nu$, i.e., $\{f(n\nu)\, |\, n\in\bz\}$. (Functions 
 $f$ with Fourier transform $\hat f$ of compact support are called {\it band-limited}. We shall adopt this convention even if $f$ is a function in several variables.)
 Shannon's theorem has now found a variety of generalizations and refinements, and it is standard fare in both Fourier analysis and in applied mathematics. Given the finite support of $\hat f$, it is well known that there is an optimal rate $\nu$ (the Nyquist rate) which gives exact reconstruction of the function from its samples. Nonetheless, there are instances where it is either desirable or unavoidable to sample beyond
the Nyquist rate, resulting in a larger, {\it oversampled} set of discrete sample points. Shannon's formula even has a formulation in Hilbert space in terms of reproducing kernels, but we shall be concerned here with wavelet bases. One way to view Shannon's interpolation is to think of the interpolation formula for the function $f$ as an expansion into a reproducing system for $L^2(\br)$, or analogously for $L^2(\br^d)$ in higher dimensions. As is well known, the reproducing systems take the form of {\it frames} (see \cite{Dau92} and \cite{Chr03}), and moreover the class of frames include wavelets, or rather wavelet bases which constitute frames, see eq (\ref{eqfr}) below. We shall be concerned with this framework for oversampling, and our results are in the context of the Hilbert space $L^2(\br^d)$. This means that our sampling points will typically constitute a rank-$d$ lattice, i.e., a discrete subgroup of $\br^d$ of rank $d$.
\section{A frame of functions on $\br^d$}
 \par
 While the early results are based on Shannon's ideas, the subject received a boost from advances in wavelets and frames, and a number of authors have recently extended and improved the classical sampling and reconstruction results. An intriguing extension is to the non-uniform case, i.e., when the sampling points are not necessarily confined to a lattice in $\br^d$. The paper \cite{AlGr01} and the book \cite{BeFe00} contain a number of such new results, and they offer excellent overviews. 
 \par
For the case of one or several real variables, there is a separate generalization of standard dyadic wavelets, again based on translation and scaling: And there is a powerful approach to the construction of wavelet bases in the Hilbert space $L^2(\br^d)$, i.e., of orthogonal bases in $L^2(\br^d)$, or just frame wavelet bases, but still in $L^2(\br^d)$. 
\par
Of course, the best known instance of this is $d = 1$, and dyadic wavelets \cite{Dau92}. In that case, the two operations on the real line $\br$ are translation by the group $\bz$ of the integers, and scaling by powers of $2$, i.e., $x\mapsto 2^j$, as $j$ runs over $\bz$.  This is the approach to wavelet theory which is based on multiresolutions analyses and filters from signal processing. In higher dimensions $d$, the scaling is by a fixed matrix, and the translations by the rank-$d$ lattice $\bz^d$. Again we will need scaling by all integral powers. We view points $x$ in $\br^d$ as column vectors, and we then consider the group of scaling transformations, $x\mapsto M^j x$  as $j$ ranges over $\bz$. 
\par
Let $M$ be a $d\times d$ dilation matrix with integer entries, such that all the eigenvalues $\lambda$ of $M$ satisfy
$|\lambda|>1$. The {\it dilation operator} induced by $M$ is $Df(x):=\sqrt{|\det M|}f(Mx)$, for $f\in\ltword$.
\par
For $u\in\br^d$, let $T_u$ denote the {\it translation operator} by $T_uf(x):=f(x-u)$.
\par
Let $H$ be a Hilbert space. A collection of vectors $\{e_i\,|\,i\in I\}$ in $H$ is called a {\it frame} if there are some constants $A,B>0$ such that, for all $f\in H$,
\begin{equation}\label{eqfr}
A\|f\|^2\leq\sum_{i\in I}|\ip{f}{e_i}|^2\leq B\|f\|^2.
\end{equation}
The constants $A$ and $B$ are called the lower and the upper frame bounds. 
\par
Our use of the term {\it oversampling} is motivated as follows: We start with a frame system in $d$ real dimensions which is based on scaling by a fixed expansive matrix $M$; and we normalize the setting such that our initial set of sample points will be located on the standard rank-$d$ lattice $\bz^d$. We
then introduce a second $d\times d$ matrix $P$ (having integer entries) and consider the larger
lattice of sample points, $P^{-1}\bz^d$. The two matrices $M$ and $P$ must satisfy a certain compatibility condition (generalizing the notion of {\it relative prime} for numbers); Definition \ref{def1_1}. In our theorem (Theorem \ref{th1}) we compare our two frame systems before and after {\it oversampling}. 
\par
We now turn to our affine frame systems, and their oversampled versions.
\par
Let $\Psi=\{\psi_1,...,\psi_n\}\subset\ltword$. The {\it affine system} generated by $\Psi$ is 
$$X(\Psi):=\{D^jT_k\psi_i\,|\,j\in\bz,k\in\bz^d,l\in\{1,...,n\}\}.$$

\par
Let $P$ be a $d\times d$ integer matrix. Denote by $p:=|\det P|$. The {\it oversampled affine system} generated by $\Psi$ relative to $P$ is
$$X^P(\Psi):=\{\frac{1}{\sqrt{p}}D^jT_{P^{-1}k}\psi_i\,|\,j\in\bz,k\in\bz^d,i\in\{1,...,n\}\}.$$

It is helpful to view sampling in context of $\ltword$-wavelets. Start with
an affine wavelet frame system in the Hilbert space $\ltword$. Our
general idea is then to represent oversampling for such an affine frame
basis in $\ltword$ as follows: While the initial frame system in $\ltword$ will have
redundancy, we show that there is a specific ``larger'' Hilbert space such
that by passing to this ambient Hilbert space, we will then get exact
(Nyquist type) sampling.  In fact we show that the ``larger'' Hilbert space
takes the form of an orthogonal sum of $\ltword$ with itself a finite number of
times $p$ say, where $p$ depends on the amount of oversampling.

\begin{definition}\label{def1_1}
The matrix $P$ is called an {\it admissible oversampling matrix} for $M$ (or simply {\it admissible}), if the matrix $PMP^{-1}$ has integer entries and $$M^{-1}\bz^d\cap P^{-1}\bz^d=\bz^d.$$
\end{definition}
\par
The second oversampling theorem states that, when $P$ is admissible, if the affine system $X(\Psi)$ is a frame for 
$\ltword$ then the oversampled affine system $X^P(\Psi)$ is a frame for $\ltword$ with the same frame bounds. This type of oversampling was introduced by Chui and Shi in \cite{CS94} for one dimension and scaling by $2$. Since then, the result has been generalized and has become known as ``the second oversampling theorem''; see \cite{RS97,CCMW02,Lau02,Jon03}. For details on the history of the second oversampling theorem we refer to \cite{Jon03}. 
\par
In this paper we will prove that, in fact, more is true: if we oversample with a matrix $P$ we will obtain frames in the larger space $H:=\underbrace{\ltword\oplus...\oplus\ltword}_{p\mbox{ times}}$. 
Specifically, we prove the following:\\
\\
{\bf Theorem \ref{th1}.} {\em  Let $d$ be given, and let two integral $d$ by $d$ matrices $M$ and $P$ satisfy the conditions in Definition \ref{def1_1}. Starting with an $M$-scale frame $X(\Psi)$ in $L^2(\br^d)$, let $X^P(\Psi)$ be the corresponding $P$-oversampled frame. Then there is an isometric and diagonal embedding of $\ltword$ in $H$ and a 
third frame in $H$ with the same frame bounds, such that $X^P(\Psi)$ arises from it as the 
projection onto the first component.  
}\\
\\

These ``super-frames'' are also generated by an affine structure on $H$. In addition, the oversampled system $X^P(\Psi)$ can be recovered as the projection of the super-frame onto the first component. Moreover, when the affine system is an orthonormal basis, the corresponding super-frame is also an orthonormal basis for $H$. 
\par
These results generalize also Theorem 5.8 in \cite{Du} which treated the case of tight frames in dimension $d=1$.

\section{Statement of the results}
\par
We will assume that $P$ is an admissible oversampling matrix for $M$. We denote by $M':=PMP^{-1}$.
\par
Let $\{\theta_r\,|\,0\leq r\leq p-1\}$ be a complete set of representatives for $P^{-1}\bz^d/\bz^d$. We can take $\theta_0=0$. The dual of this group is $(P^*)^{-1}\bz^d/\bz^d$, and let $\{\theta_q^*\,|\, 0\leq q\leq p-1\}$ be a complete set of representatives for this group. We can take $\theta_0^*=0$. The duality is given by
$$\ip{\theta_r}{\theta_q^*}=e^{2\pi iP\theta_r\cdot\theta_q^*}.$$
\par
Since $P$ is admissible, the matrix $M$ induces a permutation $\sigma$ of $\{r\,|\,0\leq r\leq p-1\}$, $M\theta_r\equiv\theta_{\sigma(r)}\mod\bz^d$ (see \cite[Proposition 2.1]{Jon03}).
\par
The dual of this map induces a permutation $\sigma^*$ of $\{q\,|\,0\leq q\leq p-1\}$, $M'^*\theta_q^*\equiv\theta_{\sigma^*(q)}\mod\bz^d$.
\par
Indeed,
$$e^{2\pi iPM\tr\cdot\ttq}=e^{2\pi iM'P\tr\cdot\ttq}=e^{2\pi iP\tr\cdot M'^*\ttq}.$$

\par
Define the Hilbert space
$$H:=\underbrace{\ltword\oplus...\oplus\ltword}_{p\mbox{ times}}.$$
On this Hilbert space we define an affine structure generated by a dilation operator $\D$, and the translation operators $\T_{P^{-1}k}$, $k\in\bz^d$.
\par
For the matrix $M$ define the unitary operators:
$$\T_{P^{-1}k}(f_q)_{0\leq q\leq p-1}=(e^{2\pi ik\cdot\ttq}T_{P^{-1}k}f_q)_{0\leq q\leq p-1},\quad(k\in\bz^d);$$
$$\D(f_q)_{0\leq q\leq p-1}=(Df_{(\sgs)^{-1}(q)})_{0\leq q\leq p-1}.$$
The operators satsify the commutation relation 
\begin{equation}\label{eqcomu0}
\D\T_{MP^{-1}k}=\T_{P^{-1}k}\D,\quad(k\in\bz^d).
\end{equation}
\par
For $f\in\ltword$, define $Sf\in H$
$$Sf=\frac{1}{\sqrt{p}}(f,...,f).$$
$S$ is an isometry.

\begin{theorem}\label{th1}
Let $\Psi:=\{\psi_1,...,\psi_n\}\subset\ltword$. The affine system $X(\Psi)$ is a frame (orthonormal basis) for $\ltword$ if and only if 
$$\mathbf{X}(S\Psi):=\{\D^j\T_{P^{-1}k}S\psi_i\,|\,j\in\bz,k\in\bz^d,i\in\{1,...,n\}\}$$
is a frame (orthonormal basis) for $H$ with the same frame bounds. The projection of $\mathbf{X}(S\Psi)$ onto the first component is the oversampled affine system $X^P(\Psi)$. 
\end{theorem}

\section{Proof of theorem \ref{th1}}

\par
We start with a lemma.
\begin{lemma}\label{lem1}
The matrix 
$$\mathcal{H}:=\frac{1}{\sqrt{p}}(e^{2\pi iP\tr\cdot\ttq})_{0\leq r,q\leq p-1}$$
is unitary.
\end{lemma}
\begin{proof}
This is a well-known fact from harmonic analysis. The matrix is the matrix of the Fourier transform on the finite abelian group $P^{-1}\bz^d/\bz^d$.
\end{proof}
\par
It will be convenient to make a change of variable $x\mapsto Px$. Let $D'$ be the corresponding dilation operator $D'f(x)=\sqrt{|\det M'|}f(M'x)$, for $f\in\ltword$.
\par
For the matrix $M'$, define the unitary operators:
$$\T'_k(f_q)_{0\leq q\leq p-1}=(e^{2\pi ik\cdot\ttq}T_kf_q)_{0\leq q\leq p-1},\quad(k\in\bz^d);$$
$$\D'(f_q)_{0\leq q\leq p-1}=(D'f_{(\sgs)^{-1}(q)})_{0\leq q\leq p-1}.$$
They satisfy the commutation relation
\begin{equation}\label{eqcomu}
\D'\T'_{M'k}=\T'_k\D',\quad(k\in\bz^d).
\end{equation}
\par

The two affine structures $\{\D,\T_{P^{-1}k}\}$ and $\{\D',\T'_k\}$ are conjugate. Define the unitary operator $U_P$ on $\ltword$ by $U_Pf(x)=\sqrt{p}f(Px)$. Then
$$U_PD'=DU_P,\quad U_PT_k=T_{P^{-1}k}U_P,\quad(k\in\bz^d).$$
\par
Define the unitary operator $\U$ on $H$ 
$$\U(f_0,...,f_{p-1})=(U_Pf_0,...,U_Pf_{p-1}).$$
Then
$$\U\D'=\D\U,\quad \U\T'_k=\T_{P^{-1}k}\U,\quad(k\in\bz^d).$$
\par
Also if 
$$S'f(x)=\frac{1}{p}(f(P^{-1}x),...,f(P^{-1}x)),\quad(f\in\ltword)$$
then $S'$ is an isometry and
$$\U S'f=Sf,\quad f\in\ltwor.$$
This shows that 
$$\U\D'^j\T'_kS'f=\D^j\T_{P^{-1}k}Sf,\quad(j\in\bz,k\in\bz^d,f\in\ltword).$$
\par
In order to show that $\mathbf{X}(S\Psi)$ is a frame with given frame bounds
it is enough to prove that 
$$\mathbf{X}(S'\Psi):=\{\D'^j\T'_kS'\psi_i\,|\,j\in\bz,k\in\bz^d,i\in\{1,...,n\}\}$$
is a frame with the same frame bounds.

\begin{lemma}\label{lem2}
Each $f\in H$ can be written uniquely as 
$$g=\sum_{r=0}^{p-1}\T'_{P\tr}S'f_r$$
for some $f_0,...,f_{p-1}\in\ltword$.
\end{lemma}
\begin{proof}
Let $g=(g_q)_q\in H$. We want 
$$(g_q(x))_q=\left(\sum_{r=0}^{p-1}e^{2\pi iP\tr\cdot\ttq}\frac{1}{p}f_r(P^{-1}(x-P\theta_r))\right)_q$$
Equivalently, in matrix form,
$$(g_q(x))_q=\mathcal{H}^t(\frac1{\sqrt{p}}f_r(P^{-1}(x-P\theta_r)))_r.$$
Using lemma \ref{lem1}, we get
$$(\frac1{\sqrt{p}}f_r(P^{-1}(x-P\theta_r)))_r=(\mathcal{H}^t)^{-1}(g_q(x))_q,$$
and this uniquely determines $f_r$, and $f_r\in\ltword$.
\end{proof}
\begin{lemma}\label{lem3}
For $f,g\in\ltword$, $k\in\bz^d$, $j\in\bz$ and $r\in\{0,...,p-1\}$, we have
$$\ip{\D'^j\T'_kS'f}{\T'_{P\tr}S'g}=\delta_{l,\sigma^j(r)}\ip{D^jT_{P^{-1}k}f}{T_{\tr}g},$$
where $l\in\{0,...,p-1\}$ is determined by $k=P\theta_l+Pm$, with $m\in\bz^d$.
\end{lemma}
\begin{proof}
First compute
$$(\D'^j\T'_kS'f)(x)=\left(\frac{1}pe^{2\pi ik\cdot\theta_{(\sigma^*)^{-j}(q)}^*}\sqrt{(\det M')^j}f(P^{-1}(M'^{j}x-k))\right)_q.$$
Therefore
$$\ip{\D'^j\T'_kS'f}{\T'_{P\tr}S'g}=
\sum_{q=0}^{p-1}\int_{\br^d}\frac{1}{p^2}e^{2\pi i(k\cdot\theta_{(\sigma^*)^{-j}(q)}^*-P\theta_r\cdot\ttq)}\cdot$$
$$\cdot\sqrt{(\det M')^j} f(P^{-1}(M'^{j}x-k))\cj g(P^{-1}(x-P\tr))\,dx.$$
But 
$$\frac{1}{p}\sum_{q=0}^{p-1}e^{2\pi i(k\cdot\theta_{(\sigma^*)^{-j}(q)}^*-P\theta_r\cdot\ttq)}=
\frac{1}{p}\sum_{q=0}^{p-1}e^{2\pi i(P\theta_l\cdot\theta_{(\sigma^*)^{-j}(q)}^*-P\theta_r\cdot\ttq)}=$$
$$\frac{1}{p}\sum_{q=0}^{p-1}e^{2\pi i(P\theta_{\sigma^{-j}(l)}\cdot\ttq-P\theta_r\cdot\ttq)}=
\delta_{\sigma^{-j}(l),r}=\delta_{l,\sigma^j(r)}.$$
We used Lemma \ref{lem1} in the second to last equality.
\par
Then 
$$\int_{\ltword} \frac{1}{p}\sqrt{(\det M')^j}f(P^{-1}(M'^{-j}x-k))\cj g(P^{-1}(x-P\tr))\,dx=$$$$
\int_{\ltword} \frac{1}{p}\sqrt{(\det M)^j}f(M^{j}P^{-1}x-P^{-1}k)\cj g(P^{-1}x-\tr)\,dx=$$$$
\int_{\ltword}\sqrt{(\det M)^j} f(M^{j}y-P^{-1}k)\cj g(y-\tr)\,dx=
\ip{D^jT_{P^{-1}k}f}{T_{\tr}g}.$$
This proves the lemma.
\end{proof}

\par
Lemma \ref{lem2} and Lemma \ref{lem3} show that the Hilbert space $H$ can be written as an orthogonal sum 
$$H=H_0\oplus H_1\oplus...\oplus H_{p-1}, \quad H_r=\T'_rS'\ltword.$$
Moreover the set $\mathbf{X}(S'\Psi)$ splits into $p$ mutually orthogonal subsets 
$$\mathbf{X}_{r}(S'\Psi):=\{\D'^j\T'_{P\theta_{\sigma^j(r)}+Pm}S'\psi_i\,|\,j\in\bz,m\in\bz^d,i\in\{1,...,n\}\}\subset H_r.$$
\par
Also note, as another consequence of Lemma \ref{lem3}, that the unitary $\D'^J$ maps $H_r$ onto $H_{\sigma^{-J}(r)}$.
\par
Therefore it is enough to show that $\mathbf{X}_r(S'\Psi)$ is a frame with the same frame bounds for $H_r$. 
We have, with Lemma \ref{lem3}, for $f\in\ltword$:
\begin{equation}\label{eqeg}
\sum_{i,j,m}|\ip{\D'^j\T'_{P\theta_{\sigma^j(r)}+Pm}S'\psi_i}{\T'_{P\tr}S'f}|^2=
\sum_{i,j,m}|\ip{D^jT_{\theta_{\sigma^j(r)}+m}\psi_i}{T_{\tr}f}|^2.
\end{equation}
But according to the proof of Theorem 3.2 in \cite{Jon03}, this quantity is bigger than 
$A\|T_{\tr}\|^2=A\|f\|^2=A\|T_{P\tr}S'f\|^2$. This yields the lower bound.
\par
For the upper bound, we have, using (\ref{eqcomu}):
$$Q:=\sum_{i,m}\sum_{j\geq0}|\ip{\D'^j\T'_{P\theta_{\sigma^j(r)}+Pm}S'\psi_i}{\T'_{P\tr}S'f}|^2=$$$$
\sum_{i,m}\sum_{j\geq0}|\ip{\D'^j\T'_{P\theta_{\sigma^j(r)-M'^jP\theta_r}+Pm}S'\psi_i}{S'f}|^2=$$$$
\sum_{i,m}\sum_{j\geq0}|\ip{\D'^j\T'_{P\theta_{\sigma^j(r)}-PM^j\tr+Pm}S'\psi_i}{S'f}|^2$$
But $P\theta_{\sigma^j(r)}-PM^j\tr=Ph$ for some $h\in\bz^d$, therefore, with Lemma \ref{lem3},
\begin{align*}
Q&=\sum_{i,m'}\sum_{j\geq0}|\ip{\D'^j\T'_{Pm'}S'\psi_i}{S'f}|^2\\
&=\sum_{i,m'}\sum_{j\geq 0}|\ip{D^jT_{m'}\psi_i}{f}|^2\leq B\|f\|^2=B\|T_{P\tr}S'f\|^2.
\end{align*}

\par
Now take $J\geq 0$. Using the fact that $\D'^J$ permutes the subspaces $H_r$, we can write 
$$\D'^J\T'_{P\tr}S'f=\T'_{P\theta_{\sigma^{-J}(r)}}S'f'$$
for some $f'\in\ltword$, with $\|f'\|=\|f\|$.
Then
$$\sum_{i,m}\sum_{j\geq-J}|\ip{\D'^j\T'_{P\theta_{\sigma^j(r)}+Pm}S'\psi_i}{\T'_{P\tr}S'f}|^2=$$$$
\sum_{i,m}\sum_{j\geq0}|\ip{\D'^j\T'_{P\theta_{\sigma^j(r)}+Pm}S'\psi_i}{\T'_{P\theta_{\sigma^{-J}(r)}}S'f'}|^2\leq B\|f'\|=B\|\T'_{P\tr}S'f\|^2.$$
Letting $J\rightarrow\infty$ we obtain the upper bound.
\par
This proves the first statement of the theorem. If $X(\Psi)$ is an orthonormal basis then $\mathbf{X}(S\Psi)$ is a frame with bounds $A=B=1$. Also the norm of the $S\Psi$ is $1$, so $\mathbf{X}(S\Psi)$ is an orthonormal basis.
\par
For the converse, take $r=0$ in (\ref{eqeg}). We know that $\mathbf{X}_0(S'\Psi)$ is a frame for $H_0$, because the families $\mathbf{X}_r(S'\Psi)$ are mutually orthogonal. Then we use (\ref{eqeg}) and conclude that $X(\Psi)$ is also a frame with the same frame bounds.
\par
The projection onto the first component $\mathcal{P}_0$ corresponds to $\theta_0^*=0$. Since $\sigma^*(\theta_0^*)=\theta_0^*$, we obtain $\mathcal{P}_0(\mathbf{X}(S\Psi))=X^P(\Psi)$. This concludes the proof of Theorem \ref{th1}.

\par
Note that from the proof of Theorem \ref{th1}, and equation (\ref{eqeg}), we obtain also the following corollary
\begin{corollary}
If $X(\psi)$ is a frame for $\ltwor$ then for all $r\in\{0,...,p-1\}$
$$\{D^jT_{\theta_{\sigma^j(r)}+m}\psi_i\,|\,j\in\bz,m\in\bz^d,i\in\{1,...,n\}\}$$ is a frame with the same bounds for $\ltword$.
\end{corollary}

\begin{remark}
 In our discussion above, we have stressed uniform sampling, i.e., the case when the sample points are located on a suitably chosen lattice in $\br^d$, and we then analyzed refinements of lattices as an instance of oversampling. But the operator theory going into our method applies also to the more general and perhaps more interesting case of non-uniform case, also called irregular sampling; see e.g., \cite{Jor83}, and \cite{BeFe00} for related modern results. Since the results are more clean in the case of lattice-refinement, we have stated our theorem in this context.
\end{remark}

\begin{acknowledgements}
The authors are pleased to thank Professor Akram Aldroubi for enlightening discussions about sampling. Both authors were supported by a grant from the National Science Foundation.
\end{acknowledgements}
\bibliographystyle{amsplain}

\end{document}